# DISTRICT LEVEL ANALYSIS OF URBANIZATION FROM RURAL-TO-URBAN MIGRATION IN THE RAJASTHAN STATE

**Dr. Jayant Singh**, Assistant Professor, Dept. of Statistics,
University of Rajasthan, Jaipur,
e-mail: jayantsingh47@rediffmail.com

**Hansraj Yadav**, Research Scholar, Dept. of Statistics,
University of Rajasthan, Jaipur,
e-mail: hansraj_yadav@rediff.com

**Dr. Florentin Smarandache**, Dept. of Mathematics,
University of New Mexico, USA,
e-mail: smarand@unm.edu

## Abstract

Migration has various dimensions; urbanization due to migration is one of them. In Rajasthan State, District level analysis of urbanization due to migrants shows trend invariably for all the districts of the state though the contribution in urbanization by the migrants varies from district to district. In some districts the share of migrants moving to urban areas is very impressive though in others it is not that much high. The migrants' contribution in urbanization is on the rising over the decades. In this paper district level migration in the Rajasthan state is examined in relation to total urbanization and urbanization due to migration.

Broadly speaking rural to urban migration is due to diverse economic opportunities across space. Throughout history migration has played substantial role in the urbanization process of several countries and still continues to play similar role. In many cases it is witnessed that more the migration higher the urbanization rate. In general, it is perceived that migration has a fairly large share in urbanization and migrants constitute a significant portion in urbanization.

At all India level rural-urban migration seems to be modest as 2001 census discloses that net rural to urban migration in 1961-71 had been 18.7 percent, in 1971-81 it was 19.6 percent, in 1981-91 migration was 21.7 percent and in 1991-01 it was 21.0 percent. So the figures reveal that there has been continuous rise in the contribution of net migration to total urban growth since the sixties though between 1991 and 2001 there has been slight decline in the rate compared to previous decade.

Migration is defined on the basis on the last residence concept hence migration in the 2001 census refers to those who migrated in ten years (1991-2001) preceding the year of survey 2001. The gross decadal inflow of rural to urban migrants as a percentage of total urban population in 2001 turns out to be a little above 7 per cent at the all-India level (Table on next page). However, it varies considerably across states. Both industrialized states like Gujarat and Maharashtra and the backward states like Orissa and Madhya Pradesh show high rates of migration. Similarly examples can be found from both the types of states which have recorded sluggish migration rate, e.g. industrialized states such as Tamilnadu and West Bengal and backward states such as Uttar Pradesh, Bihar and Rajasthan. This reveals that share of the rural to urban migrants in urbanization differs from state to state. A table giving rural to urban migrants for the period 1991-2001 as a % of urban population relation is given in table on ensuing page.

**Table 1**
**Rural-Urban migration for 1991-2001 as a % of urban population**

| States | Rural-to-Urban Migrants (1991-2001) as a % of Urban Population |
|---|---|
| Andhra Pradesh | 6.72 |
| Assam | 7.12 |
| Bihar | 6.28 |
| Gujarat | 10.63 |
| Haryana | 11.45 |
| Karnataka | 7.03 |
| Kerala | 6.99 |
| Madhya Pradesh | 9.50 |
| Maharashtra | 10.41 |
| Orissa | 10.97 |
| Punjab | 7.63 |
| Rajasthan | 6.18 |
| Tamil Nadu | 3.34 |
| Utter Pradesh | 4.44 |
| West Bengal | 4.83 |
| All India | 7.32 |

Source: Census of India 2001, Migration Tables.

Nevertheless, rural-urban migration rates at intrastate level have been a phenomenal in India as this flow dominates the interstate flows. Since the intrastate migration rates are much higher in magnitude than the interstate migration rates therefore it makes an interesting subject area to comprehend various economic, social and cultural factors connected closely with it. A district level analysis for Rajasthan state is thus attempted to perceive urbanization due to migration their interlinkages and affiliations.

**Urbanization Trend in the State of Rajasthan**

According to the census report of 2001 the share of urban population in Rajasthan has inched up to 23.38% as compared to 15.06% mentioned in the census report 1901. Number of towns in the Rajasthan increased to 216 in the census 2001 against 133 in the 1901 census that depicts 62.4% of growth in this period of time whereas at national level this growth has been 169.36% in same time span. Share of Rajasthan's urban population in the country dropped to 4.6% from 5.98% over a period of century whereas in terms of growth of number of towns, state share also slipped down to 4.18% from 6.94% in this same period of time. Therefore, it can be clearly claimed that Rajasthan has to go a long way to match with national figures as regards the characteristics of urbanization is concerned whether it is growth in urban population or towns. However, there has been a meager improvement in the percentage share of state's urban population in the national urban population as it has grown to 4.1% to 4.52%, 4.52% to 4.62% and then to 4.64% in last three successive censuses.

**District Level Analysis for Rajasthan**

The migrants contribution in urbanization is on the rising over the decades as 16.4% of the total migrants in the Rajasthan settled in urban areas during the period 1971-80 and the figure which went up to 22.4% during the duration 1981-1990 and further advanced to 25.4% in the duration 1991-2000.  This trend is evident invariably in all the districts of the state though the contribution in urbanization by the migrants varies from district to district. In some districts the share of migrants moving to urban areas is very impressive though in others it is not that much high.

The census analysis of Barmer district of Rajasthan reveals that 7.7%, 7.1% & 4.0% of total migrants moved to urban areas in last three decades i.e. 1991-2000, 1981-90 & 1971-1980. This percentage share for Jalore was 9.6, 8.1 & 4.7%, and for Banswara it was 9.1, 7.9 & 4.7%.  The figures disclose that these districts had poor share of migrants to urban areas. On the other hand there are districts like Jaipur, Ajmer, Kota and Bhilwara where the percentage share of migrants settling in urban areas with context to the total migrants is comparatively much higher. This percentage share of rural migrants in three last successive decades for these districts is given in table placed below.

**Table 2**

**Share of Rural Migrants in selected Districts during last three decadal period**

| District / period | 1991-2000 | 1981-90 | 1971-1980 |
|---|---|---|---|
| Kota | 56.8 | 54.3 | 50.7 |
| Jaipur | 53.2 | 48.5 | 35 |
| Ajmer | 41.4 | 35.6 | 28.7 |
| Bhilwara | 31.1 | 25.0 | 14.8 |
| Jodhpur | 26.8 | 18.7 | 12.4 |

To apprehend the trends in the migration of population to the urban areas in different districts of Rajasthan, based on the share of urbanization due to migration can be categorized as follows:

Category 1: Higher During all the three decades

Category 2: Higher during 1991-2000 & 1981-91 but lower in 1971-80

Category 3: Higher during 1991-2001 but lower in last two decades

Category 4: Lower During all the three decades

Category 5: Lower during 1991-2000 & 1981-90 but higher in 1971-80

Category 6: Lower during 1991-2000 but higher in 1981-90 & 1971-80

Districts falling in Category 1 are those, which observed higher urbanization due to migration in comparison with state level figures during three consecutive decadal periods. In these Districts, the proportion of migrants coming to urban areas is higher than the state proportion of such migration.

District falling in Category 2 performed better as far urbanization due to migration in last two decades is concerned. Districts in this category observed higher urbanization share due to migration than what was seen in the state in last two decadal times whereas three decades back share of migrants to urban areas was lower in these districts from that of state in overall.

Similarly, Category 3 is featuring districts that have observed higher urbanization share due to migration than to state in recent decade though that particular district was falling below than state share in two previous consecutive decades.

Category 4 to 6 are counterpart of category 1 to 3 where share of migrants moving to urban areas in total migrants for a district is lower than state share of migrants moving to urban areas as regards total migrants of the state.

**Table 3**

**Classification of District according to Urbanization Trends in last three decades**

| Category | Districts |
|---|---|
| Category 1 | Ganganagar, Bharatpur, Swaimadhopur, Jaipur, Pali, Ajmer, Kota |

| Category 2 | Bhilwara |
|---|---|
| Category 3 | Jodhpur |
| Category 4 | Alwar, Dholpur, Karauli, Dausa, Sikar, Nagaur, Barmer, Jalore, Sirohi, Tonk, Bundi, Rajsamand, Udaipur, Dungarpur, Banswara, Baran, Chittorgarh, Jhalawar |
| Category 5 | Bikaner, Jhunjhunu |
| Category 6 | Hanumangarh, Churu |

Classification elucidated above undoubtedly depicts that there are only seven districts where there is larger urbanization due to rural migrants in context with the overall state level migration and urbanization figures over three consecutive decades. Notwithstanding there are 18 districts having lower urbanization due to migration than to state level migrant urbanization.

2001 census report explains that Jodhpur is the only district where urbanization due to migration has improved with regard to the figure of state in total. Similarly, district Bhilwara has witnessed this edge in two recent decades. In two recent decades there is improvement in the data of Bhilwara in relation to urbanization due to migration is concerned decades otherwise three decades back the urbanization due to migration for Bhilwara was lower than state figures. Jodhpur showed this improvement in last decade even though it was lagging behind in two previous decades.

Bikaner and Jhunjunu are way behind in showing any improvement in urban migration to state share in last two decades while Hanumagarh & Churu showed no improvement only in last decade. Jaiselmer is the district that doesn't observe any clear-cut pattern on account of migrants share in relation to state.

**Urbanization and Migration:**

It is well evident that number of rural migrants as regards total migrants is considered as an extent of urbanization by migration in a particular category. Districts are classified in the groups where percentage of migrants attributing to urbanization is <20%, 20-50 and >50% in all the three durations 1971-80, 1981-90 and 1991-2000 and the result is summarized as below:

**Table 4**

**Number of Districts according to range of Urbanization in last three Census**

| Range of urbanization (in%) | 2001 | 1991 | 1981 |
|---|---|---|---|
| | Number of Districts | | |
| <20 | 10 | 16 | 28 |
| 20-50 | 20 | 14 | 3 |
| >50 | 2 | 2 | 1 |

Its is evident from above classification that there is stark variation in the urbanization by migrants in various census barring the category of the districts that are having >50% of urban migrants in total migrants as there are only district since last two census against one three decade back where as considerably shift in the other two categories of 20-50% and <20% urbanization due to migration is there in this three decadal period. There are more districts classified in the category 20-50% during the recent decades whereas the number of districts in the category <20% has gone down in the recent decades.

**Comparative Analysis of Total Urbanization & Urbanization due to Migration:**

Migration is an important part of the urbanization and in many cases it is attributing predominately in the urbanization. Indicator of rate of Urbanization can be defined as below:

1. Total Urbanization rate: is the percentage of population living in urban areas to the total population

2. Urbanization rate due to migration: is the percentage share of rural migrants to the total migrants.

The result of the comparative investigation made on the basis of above mentioned two indicators for the last decadal period i.e. 1991-2001 is examined in coming paragraphs.

State urbanization rate is the share of urban population to the total population at state level and similarly it is counted on districts level. Consequently these two rates are compared at state and districts level to analyze the urbanization trend and to establish its association with the migration. At state level 23.4% of the total population is urbanized and 22.9% of migrants are coming to urban areas thus at state level the urbanization rate through migrants is compatible to the total urbanization rate. Barmer and Jalore are two districts in which urbanization through migrants' rate is below 20% as the urbanization rate of the migrants to these districts is mere 15 & 19%.

Rate of urbanization through migrants in Jaipur is (73.6%), Kota (68.2%), Ajmer (53.8%) and Udaipur (50%) and thus these districts have more than 50% of rural migrants and this can be summed up as more than half of the migrants to these districts are settling in urban areas. Bikaner and Churu are the only districts observed where urbanization through migrants rate is lower than total urbanization rate of the state. This difference was more than 32% for the Udaipur and Banswara districts and for seven districts it was more than 20%. The classification of number of districts based on the range of these two urbanization indicators is classified in coming table.

**Table 5**

**Total Urbanization Rate vis-à-vis Urbanization Rate due to Migration**

| Range of Urbanization rate | | >50% | 40-50% | 30-40% | 20-30% | <20% |
|---|---|---|---|---|---|---|
| Combined (Male & female) | **Total Urbanization rate** | 1 | 2 | 2 | 8 | 19 |
| Male | | 1 | 1 | 2 | 9 | 19 |
| Female | | 1 | 1 | 3 | 7 | 20 |
| Combined (Male & female) | **Urbanization rate due to migration** | 4 | 5 | 8 | 13 | 2 |
| Male | | 12 | 8 | 4 | 9 | 12 |
| Female | | 2 | 2 | 11 | 10 | 7 |

Clearly the migration witnesses a better urbanization rate and there are more districts classified in higher range of urbanization rates than the number of district classified in lower range in accordance with the total urbanization rate of the districts.

Technique of non-parametric test is used for district level analysis of the urbanization to examine the migration to different districts having same size of population. District are ranked on the basis of the total urban population and urban population due to migration and these formed two groups of Non-parametric test and Wilcoxon - Mann/Whitney Non parametric Test is employed for equality of K universes for total population and Male & Female population and results of the analysis done in Megastat is as below:

| TOTAL | | |
|---|---|---|
| n | sum of ranks | |
| 32.00 | 698.00 | Group 1 |
| 32.00 | 1382.00 | Group 2 |
| 64.00 | 2080.00 | Total |
| | 1040.00 | expected value |
| | 74.48 | standard deviation |
| | -4.59 | Z |

|   |         | p-value (two-tailed) |
|---|---------|----------------------|
| MALE | | |
| n | sum of ranks | |
| 32.00 | 612.00 | Group 1 |
| 32.00 | 1468.00 | Group 2 |
| 64.00 | 2080.00 | Total |
|  | 1040.00 | expected value |
|  | 74.48 | standard deviation |
|  | -5.74 | Z |
|  | 0.00 | p-value (two-tailed) |
| FEMALE | | |
| n | sum of ranks | |
| 32.00 | 775.00 | Group 1 |
| 32.00 | 1305.00 | Group 2 |
| 64.00 | 2080.00 | Total |
|  | 1040.00 | expected value |
|  | 74.48 | standard deviation |
|  | -3.55 | Z |
|  | .0004 | p-value (two-tailed) |
| GROUP 1 URBANIZATION IN TOTAL POPULATION | | |
| GROUP 2 URBANIZATION BY MIGRATION | | |

Clearly, the above examined district level analysis reveals that total urbanization and urbanization due to migration differs significantly. Male and female population and districts have significant impact on total urbanization & urbanization due to migration. Thus the relative magnitude of total urbanization and urbanization due to migration differs significantly for the districts for both genders and combined.